# EXPLORING THE RELATION BETWEEN MATHEMATICAL VALUES AND ACHIEVEMENT AMONG GIRLS: A COMPARATIVE ANALYSIS IN SINGLE-SEX VS. CO-EDUCATIONAL SETTINGS USING TIMSS 2019 NZ DATA

Huayu Gao, Tanya Evans & Gavin T. L. Brown

University of Auckland, New Zealand

*Grounded in the Social Cognitive Career Theory, this study investigates the influence of values on girls' mathematics achievement across socio-economic status (SES) settings, contrasting single-sex and co-educational schools. ANOVA on 2019 TIMSS New Zealand data (n = 2,898) reveals that single-sex education is associated with enhanced girls' mathematical values in high-SES settings. However, its effect on translating these values into improved mathematics performance is relatively limited. Affluent learning resources are the more immediate factor in improving mathematics performance. Moreover, in low-SES environments, the relation between values and mathematics achievement exhibits a complex nonlinear pattern.*

## INTRODUCTION

In most OECD countries, including New Zealand, the proportion of female graduates in traditionally male-dominated fields of Science, Technology, Engineering, and Mathematics (STEM) in higher education is significantly lower than that of males, with an average of only 33% (OECD, 2023). This gender segregation trend extends to the occupational domain. For instance, a New Zealand government report indicates that only 27% of digital technology positions are held by women (Hindle & Muller, 2021). Academia faces similar challenges; according to data from the University of Auckland (NZ), the distribution of female faculty ranks in the Faculty of Science and Engineering is significantly lower than that of males, particularly in engineering, where the gap is most prominent (Brower & James, 2020). Without effective interventions, this gender imbalance is projected to persist until 2070 (Brower & James, 2020).

Girls in some Islamic countries exhibit a notable advantage in math achievement, values, and career aspirations for math-intensive careers, which stands in stark contrast to the gender difference trends observed in most Western regions (Michaelides et al., 2019). This phenomenon may be attributed to the unique social circumstances faced by women in these Islamic countries. Women's social status is relatively low in many Islamic nations, while careers in STEM fields often offer higher salaries. Consequently, girls may aspire to excel in mathematics and subsequently pursue high-paying math-related careers as a means to elevate their social standing (Stoet et al., 2018). This observation suggests that among the female population, there may be a significant association between value, mathematics performance, and the likelihood of entering math-intensive fields. From this perspective, girls' values and achievements could serve as ideal entry points for investigating their propensity to pursue math-related careers.

Apart from social factors, the single-sex education widely adopted in these Islamic countries—where girls attend all-girls schools taught by female teachers and boys attend all-boys schools taught by male teachers (Michaelides et al., 2019)—may also contribute to girls' mathematics achievement to a certain extent. However, when examining the impact of single-sex education on girls' mathematics performance, it is crucial not to overlook the significant role of students' socioeconomic status (SES). SES is considered the most influential factor affecting students' mathematics achievement, accounting for more than 50% of the variance in girls' mathematics performance (Gao et al., 2025). Relevant research indicates that when SES is not controlled for, gender differences in mathematics achievement exist; however, when SES is controlled for, these gender differences significantly narrow or disappear (Clavel & Flannery, 2023). Therefore, to more accurately assess the effectiveness of single-sex education, it is necessary to consider and control for the SES. Smith and Evans' (2023) study in New Zealand serves as an excellent example. After controlling for students' SES, they found that girls from low SES families in single-sex schools outperformed girls attending co-educational schools in mathematics ($g = 0.90$), further supporting the potential advantages of single-sex education.

New Zealand's education system offers a distinctive setting for investigating the effects of single-sex schools. Among the 374 secondary schools in NZ, approximately 16% are single-sex girls' schools, with a striking 91% of these institutions being state or state-integrated schools (MoE, 2023). This contrasts sharply with the United States and Australia, where single-sex schooling is primarily concentrated in private or Catholic schools. Furthermore, New Zealand employs a decile system to measure the socioeconomic background of schools. This system ranks schools from 1 to 10, with decile 1 representing schools whose students come from the lowest 10% of families and decile 10 representing schools whose students come from the highest 10% of families.

**THEORETICAL BACKGROUND**

The performance model of the Social Cognitive Career Theory (SCCT) provides a comprehensive framework for understanding the interplay between individual motivation and environmental factors in shaping math-related performance (Lent & Brown, 1994). More specifically, the model centers on the core constructs of self-efficacy and outcome expectation, which dynamically interact with person factors, contextual influences, and learning experiences to collectively shape academic outcomes. Outcome expectations, as a center component, encompass the concept of values, reflecting individuals' preferences for the "reinforcers" of academic activities (Dawis & Lofquist, 1984). For instance, mathematical values are not only embodied in the interest in mathematics itself but also include the expectation that learning mathematics can lead to career opportunities, social recognition, and economic rewards. Contextual and experiential factors act as antecedents to the formation of values or directly facilitate or hinder individuals' mathematics performance, with gender serving as a significant moderator in these relationships. For example, girls may encounter a higher prevalence of gender role stereotypes and negative evaluations, such as the pervasive notion that girls lack aptitude in mathematics, which can erode their outcome expectations and, consequently, impinge upon their mathematics performance.

## RESEARCH QUESTIONS

This study addresses the relationship between school SES (e.g., decile), school gender (e.g., single-sex vs coeducational), mathematical values (e.g., self-reported values of mathematics), and mathematics achievement in TIMSS for Year 9 girls in New Zealand. Two questions were formulated:

1. Is there a disparity in girls' mathematical values between single-sex and co-ed schools?

2. What is the relationship between mathematics values and mathematics achievement, contrasting girls in single-sex and co-ed schools?

## METHODS

The study comprised 2898 Year 9 girls in the first year of New Zealand secondary education from the 2019 TIMSS dataset. The 2023 TIMSS dataset was excluded due to its limitations: the impact of COVID-19-related disruptions, such as online learning, complicates disentangling the effects of school factors on mathematics achievement. Furthermore, access to key non-public variables, such as school gender and decile, was unavailable at the time of this research, further reducing the dataset's applicability for this study.

Decile and school gender are critical variables in this study. Decile bands are low = deciles 1-3, medium = deciles 4-7; high = deciles 8-10. Of these participants, 906 were enrolled in single-sex girls' schools, with 9.2%, 35.8%, and 55.1% across low-decile, medium-decile, and high-decile schools, respectively. Additionally, 1992 girls attended co-ed schools, distributed as 26.0%, 41.6%, and 32.5% across low-decile, medium-decile, and high-decile schools. Based on percentages, the distribution clearly indicates a substantial advantage of single-sex schools in terms of SES ($\chi^2 = 14.21$, $p < .001$).

The data were subjected to two-way ANOVA and three-way ANOVA to investigate the relationships between decile band (DB: Low, Medium, and High), school gender (SG: Single-sex and Co-ed schools), mathematical values (MV: Low, Medium, and High), and girls' mathematics achievement (MA).

Given the substantial disparity in girls' enrollment between co-ed and single-sex schools, potential bias in ANOVA was overcome by drawing 20 samples, mirroring the enrollment figures of single-sex schools stratified by DB. Each sample comprised 232, 376, and 294 girls from low, medium, and high decile schools. The null hypothesis is rejected if only one of the 20 iterations proves statistically insignificant.

**RESULTS**

A two-way ANOVA was conducted to investigate the impact of school gender and decile settings on values in mathematics. A statistically significant interaction was observed between SG and DB on MV (Table1: $F_{(4, 2858)} = 5.523$, $p = .004$, $\eta p2 = .004$).

| Source | df | F | Sig. | ηp2 |
|---|---|---|---|---|
| Corrected Model | 5 | 3.113 | 0.008 | 0.005 |
| Intercept | 1 | 42743.54 | <.001 | 0.937 |
| SG | 1 | 0.061 | 0.804 | 0.000 |
| DB | 2 | 2.975 | 0.051 | 0.002 |
| SG * DB | 2 | 5.523 | 0.004 | 0.004 |

Table 1: Interaction between gendered school and decile setting

Consequently, further analysis of simple main effects for SG and DB was performed, employing Bonferroni-adjusted statistical significance at the $p < .025$ level (Table 2). In terms of SG, a significant difference in mean "MV" scores for single-sex schools was found ($F_{(2, 2870)} = 4.471$, $p = .012$). A significant difference in mean "MV" scores between single-sex and co-ed schools for low-decile was observed ($F_{(1, 2870)} = 6.048$, $p = .014$). Among the 20 samples, only one was deemed insignificant, with the remaining 19 demonstrating statistical significance at $p < 0.05$.

| DB | df | F | Sig. | SG | df | F | Sig. |
|---|---|---|---|---|---|---|---|
| Low-DB | 1 | 6.048 | 0.014 | Single | 2 | 4.471 | 0.012 |
| Medium-DB | 1 | 3.196 | 0.074 | Co-Ed | 2 | 1.908 | 0.149 |
| High-DB | 1 | 5.221 | 0.022 | | | | |

Table 2: Main effect of decile setting and gendered school

Pairwise comparisons were conducted within each main effect, with Bonferroni-adjusted confidence intervals and p-values. Girls scored significantly lower mean "MV" in single-sex than co-ed schools in low-decile settings (0.540; 95% CI [0.110,0.971]; $p = .014$). Furthermore, girls had significantly lower mean "MV" scores in single-sex schools in low-decile compared to those in medium and high-decile bands.

In conclusion, a notable disparity in values between single-sex and co-ed schools was evident in low-decile bands, with the latter exhibiting significantly higher values. Moreover, across the low-to-high decile spectrum, the values of girls in single-sex schools increased, while there was no significant difference in values for girls in co-ed schools.

A three-way ANOVA examined the nuanced impacts of school gender, decile band, and values on mathematics achievement. A statistically significant three-way interaction among these factors was identified ($F_{(4, 2858)}$ = 2.398, $p$ = .048), with a standardized effect size ($r$) of 0.064 (95% CI [0.0494, 0.0787]). Notably, Cohen's (1988) guideline designates this effect size as consistently ignorable (Table 3).

|  | df | F | Sig. | $\eta p^2$ |
|---|---|---|---|---|
| Intercept | 1 | 52885.17 | <.001 | 0.949 |
| SG | 1 | 39.557 | <.001 | 0.014 |
| DB | 2 | 84.87 | <.001 | 0.056 |
| MV | 2 | 21.201 | <.001 | 0.015 |
| SG * MV | 2 | 1.343 | 0.261 | 0.001 |
| SG * DB | 2 | 11.604 | <.001 | 0.008 |
| DB * MV | 4 | 4.931 | <.001 | 0.007 |
| SG * DB * MV | 4 | 2.398 | 0.048 | 0.003 |

Table 3: Interaction between gendered school, decile setting, and values

Significant two-way interactions were found between SG and DB ($F_{(2, 2858)}$ = 11.604, $p$ < .001) and between MV and DB ($F_{(4, 2858)}$ = 4.931, $p$ < .001) (Table 4). Furthermore, the main effects of SG existed in low-decile bands and MV in high-decile bands. Pairwise comparisons showed that girls from single-sex schools had significantly higher mean "MA" scores than their co-ed counterparts in low-decile bands (61.659; 95% CI [41.713, 81.606]). High MV girls demonstrated markedly higher mean "MA" scores in high decile band compared to their counterparts with low and medium MV (Δlow: 64.206; 95% CI [48.008, 80.403] and (Δmedium: 35.446; 95% CI [23.187, 47.706].

|  | SG | | | | MV | | | |
|---|---|---|---|---|---|---|---|---|
|  | df | F | Sig. | $\eta p^2$ | df | F | Sig. | $\eta p^2$ |
| Low-DB | 1 | 36.739 | <.001 | 0.013 | 2 | 1.43 | 0.239 | 0.001 |
| Mid-DB | 1 | 2.479 | 0.116 | 0.001 | 2 | 8.576 | <.001 | 0.006 |
| High-DB | 1 | 3.524 | 0.061 | 0.001 | 2 | 49.482 | <.001 | 0.033 |

Table 4: Main effect of gendered school and values

Consequently, a significant disparity in girls' mathematics achievement emerged between single-sex and co-ed schools in low-decile bands. Notably, values did not significantly influence girls' math scores in low-decile and medium-decile bands. However, a positive relationship emerged in high-decile bands, indicating that higher MV corresponded to higher MA only among those girls.

## DISCUSSION

### Dual Necessities for Enhancing Girls' Mathematical Value: Resource Guarantee and Single-sex Education

As the decile level increases, the development of girls' values in co-educational schools remains relatively stagnant, whereas in single-sex schools, these values exhibit significant growth. This finding indicates that cultivating motivation, such as mathematical values, necessitates the availability of abundant resources as a foundational prerequisite. Within this context, single-sex education can effectively maximize its potential impact. The increase in decile levels implies that families and schools can provide more abundant resources, such as high-quality teaching equipment, diverse curriculum choices, and more extracurricular activity opportunities. Single-sex schools, compared to co-ed schools, can more effectively transform these resources into personalized support, helping girls break gender stereotype roles, such as "mathematics is not for girls," thereby promoting the development of girls' values (Smith & Evans, 2023). The initiatives of St Cuthbert's College for girls (decile 10) in New Zealand serve as a prime example. The college regularly invites women who have achieved outstanding accomplishments in STEM fields as role models to challenge prevalent gender stereotypes, thereby enhancing girls' perceptions of mathematics. In contrast, co-ed schools may struggle to provide targeted support and activities for girls. Therefore, the long-standing gender role expectations within co-ed schools are challenging to eliminate, weakening the promoting effect of affluent resources on girls' values in the co-ed environment.

### Translating Mathematical Values into Achievements: The Significance of Learning Resources and the Limited Role of Single-Sex Education

In high decile schools, the positive impact of values on girls' mathematics achievement is prominently manifested. However, this positive impact does not appear to exhibit significant differences between single-sex and co-ed schools. This phenomenon suggests that the process of transforming motivation into academic achievement needs to be supported by sufficient resources, while gender-related educational influences are relatively limited. Nevertheless, it would be one-sided to completely deny the impact of single-sex education on math achievement based solely on this phenomenon. As a fundamental condition for educational quality, sufficient resources have a more direct and immediate effect on improving mathematics performance and can even explain 56% of the variance in mathematics achievement (Gao et al., 2025). Moreover, through the above analysis, it can be found that single-sex education can effectively promote the enhancement of girls' mathematical values, while the promoting effect of motivation on mathematics achievement has a certain lagging effect. The immediate effect of resources and the lagging effect of motivation also explain the research results obtained by Smith and Evans (2023) in cross-sectional data: when both single-sex and co-educational schools have abundant resources, the performance differences of girls in mathematics and science are not significant. Therefore, considering the unique function of single-sex education in promoting girls' values, we cannot completely deny the role of single-sex education in girls' mathematics achievement. Further longitudinal analyses are imperative to elucidate this issue and provide a more comprehensive understanding of the intricate interplay between educational settings, values, and academic performance.

**Complex Relations in Low-decile Band: The Nonlinear Influence between Mathematical Values and Achievement**

In the low-decile band, the advantages of single-sex education in fostering girls' mathematical values appear to be diminished. In economically disadvantaged settings, the potential benefits of single-sex education in fostering values are inevitably overshadowed by systemic issues arising from resource scarcity. Low-decile schools often face challenges such as inadequate facilities, limited teaching resources, and weak socio-economic support networks. In such contexts, girls may not receive sufficient additional support and compensatory education. The targeted attention and support that single-sex education aims to provide may be weakened by the constraints imposed by resource scarcity, consequently affecting girls' development of values in mathematics.

Regarding math achievement, girls in single-sex schools surpass their couterparts in co-ed settings. Under resource-constrained conditions, single-sex schools appear to be more efficacious in enhancing girls' math performance. However, this does not necessarily imply that gender-specific support activities have a significant impact on math achievement. In fact, optimized resource allocation in single-sex schools may be the more direct factor influencing girls' math performance. In co-ed schools, girls' time and resources for math learning are often compromised due to boys' higher frequency of interaction with teachers. This imbalanced resource distribution undermines girls' mathematical confidence and engagement (Smith & Evans, 2024). Furthermore, the discrepancy between higher mathematical values and lower confidence in abilities in co-ed settings may elicit negative emotions, such as anxiety and helplessness (Pekrun, 2006). These emotions further exacerbate the detrimental effects of co-ed environments on math performance. In contrast, single-sex schools, by virtue of the absence of boys, can concentrate learning resources on girls. The direct impetus of these learning resource allocations on math achievement can effectively mitigate the negative impacts on academic performance caused by resource scarcity.

Notably, the relationship between values and mathematics performance exhibits complex non-linear patterns in low-decile environments. In other words, no significant correlation is observed between values and math performance in low-SES contexts. This suggests that economic disadvantage constrains the effective translation of mathematical value into achievement, challenging numerous academic success theories, including SCCT. These theories typically assume the universal applicability of motivation to academic achievement, but their explanatory power may be primarily limited to academic contexts with resource-rich or high-achieving groups.

**REFERENCES**

Brower, A., & James, A. (2020). Research performance and age explain less than half of the gender pay gap in New Zealand universities. *Plos One*, *15*(1), e0226392. https://doi.org/10.1371/journal.pone.0226392

Clavel, J. G., & Flannery, D. (2023). Single-sex schooling, gender and educational performance: Evidence using PISA data. *British Educational Research Journal*, *49*(2), 248–265. https://doi.org/10.1002/berj.3841

Cohen, J. (1988). Set correlation and contingency tables. *Applied Psychological Measurement*, *12*(4), 425–434. https://doi.org/10.1177/014662168801200410

Dawis, R. V., & Lofquist, L. H. (1984). *A psychological theory of work adjustment: an individual-differences model and its applications*. University of Minnesota Press.

Gao, H., Evans, T., & Brown, G. T. L. (2025). Exploring gender differences in tertiary mathematics-intensive fields: A critical review of social cognitive career theory. *Proceedings of the 27th Annual Conference on Research in Undergraduate Mathematics Education*. Alexandria, VA.

Hindle, S., & Muller, G. (2021). *Digital skills Aotearoa: Digital skills for our digital future*. MBIE.

Lent, R. W., Brown, S. D., & Hackett, G. (1994). Toward a unifying social cognitive theory of career and academic interest, choice, and performance. *Journal of Vocational Behavior*, *45*(1), 79–122. https://doi.org/10.1006/jvbe.1994.1027

Michaelides, M., Brown, G. T. L., Eklöf, H., & Papanastasiou, E. (2019). *Motivational profiles in TIMSS mathematics: Exploring student clusters across countries and time*. Springer Open & IEA. https://doi.org/10.1007/978-3-030-26183-2

MoE. (2023). *Education Counts*. Retrieved from https://www.educationcounts.govt.nz/statistics/school-rolls

OECD. (2023). *Education at a Glance 2023 Sources, Methodologies and Technical Notes*. OECD Publishing, Paris. https://doi.org/10.1787/d7f76adc-en

Pekrun, R. (2006). The control-value theory of achievement emotions: Assumptions, corollaries, and implications for educational research and practice. *Educational Psychology Review*, *18*, 315–341. https://doi.org/10.1007/s10648-006-9029-9

Smith, A., & Evans, T. (2024). Gender gap in STEM pathways: The role of gender-segregated schooling in mathematics and science performance. *New Zealand Journal of Educational Studies*, 1–19. https://doi.org/10.1007/s40841-024-00320-y

Stoet, G., & Geary, D. C. (2018). The gender-equality paradox in science, technology, engineering, and mathematics education. *Psychological Science*, *29*(4), 581–593. https://doi.org/10.1177/0956797617741719